# DEFORMATION and EXTENSION of FIBRATIONS of SPHERES by GREAT CIRCLES

## Patricia Cahn, Herman Gluck and Haggai Nuchi

In a 1983 paper with Frank Warner, we proved that the space of all great circle fibrations of the 3-sphere $S^3$ deformation retracts to the subspace of Hopf fibrations, and so has the homotopy type of a pair of disjoint two-spheres. Since that time, no generalization of this result to higher dimensions has been found, and so we narrow our sights here and show that in an infinitesimal sense explained below, the space of all smooth oriented great circle fibrations of the 2n+1 sphere $S^{2n+1}$ deformation retracts to its subspace of Hopf fibrations. The tools gathered to prove this also serve to show that every germ of a smooth great circle fibration of $S^{2n+1}$ extends to such a fibration of all of $S^{2n+1}$, a result previously known only for $S^3$.

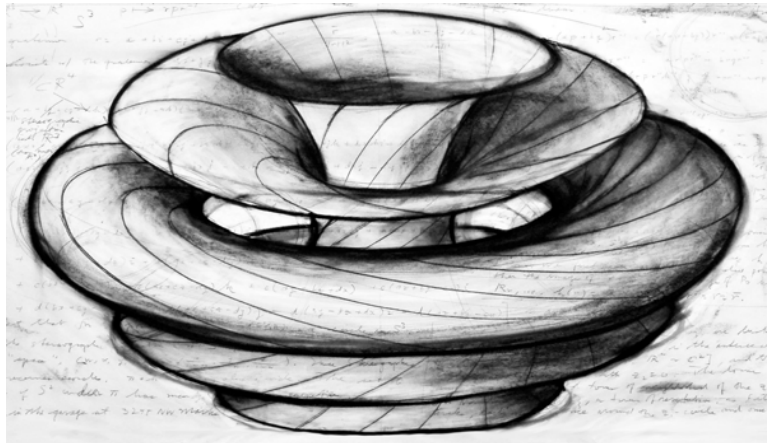

**Hopf fibration of 3-sphere by great circles**
**Lun-Yi Tsai    Charcoal and graphite on paper    2007**

# INTRODUCTION

Consider a fibration  F  of the unit 2n+1 sphere  $S^{2n+1}$  by oriented great circles, and focus on one of the fibres  P ,  as shown below.

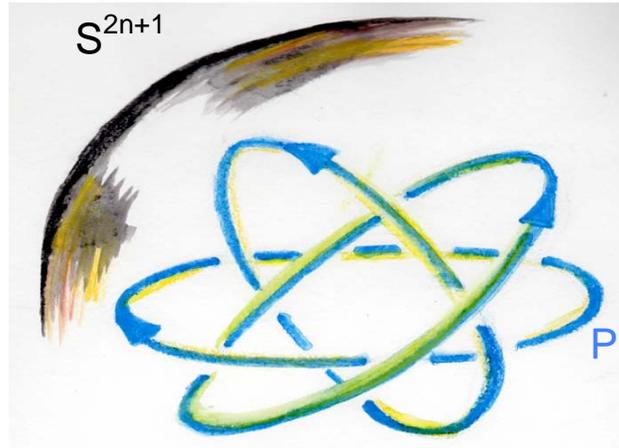

**Figure 1.  A fibration of  $S^{2n+1}$  by oriented great circles**

The oriented great circle  P  spans an oriented 2-plane through the origin in  $R^{2n+2}$ ,  which we also denote by  P ,  and so appears as a single point in the Grassmann manifold  $G_2R^{2n+2}$  of all such oriented 2-planes.  If the fibration  F  is smooth, then its base space  $M_F$  appears as a smooth 2n-dimensional submanifold of this Grassmann manifold, and we can focus on the tangent 2n-plane  $T_PM_F$  to  $M_F$  at  P .

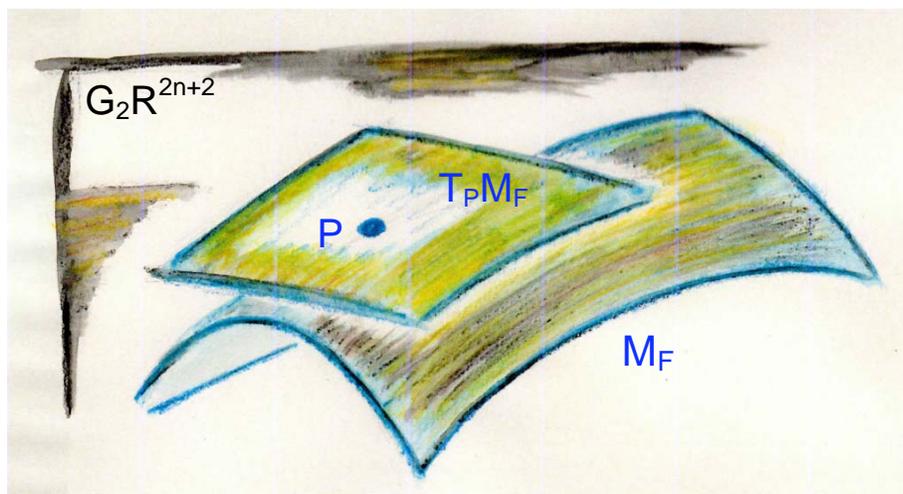

**Figure 2.  The base space  $M_F$  of the fibration  F ,
and its tangent plane  $T_PM_F$  at  P .**

**THEOREM A.** The space $\{T_P M_F\}$ of tangent 2n-planes at $P$ to the base spaces $M_F$ of all smooth oriented great circle fibrations $F$ of $S^{2n+1}$ containing $P$, deformation retracts to its subspace $\{T_P M_H\}$ of tangent 2n-planes to such Hopf fibrations $H$ of $S^{2n+1}$.

**THEOREM B.** Every germ of a smooth fibration of $S^{2n+1}$ by oriented great circles extends to such a fibration of all of $S^{2n+1}$.

A **germ** of a fibration of $S^{2n+1}$ by oriented great circles consists of such a fibration in an open neighborhood of a given fibre $P$, with two germs equivalent if they agree on some smaller neighborhood of $P$. To **extend** such a germ to a fibration of all of $S^{2n+1}$ means to find a fibration of all of $S^{2n+1}$ which agrees with the given germ on some neighborhood of $P$.

The path to the above theorems consists of the following steps.

First two definitions. The **bad set** $BS(P) \subset G_2 R^{2n+2}$ consists of all oriented 2-planes through the origin in $R^{2n+2}$ which meet $P$ in at least a line, and the **bad cone** $BC(P) \subset T_P(G_2 R^{2n+2})$ is its tangent cone at $P$.

**PROPOSITION 1.** A closed smooth 2n-dimensional submanifold $M$ of $G_2 R^{2n+2}$ is the base space of a smooth fibration of $S^{2n+1}$ by great circles if and only if it is transverse to the bad cone at each of its points.

Next we focus in on the tangent space $T_P(G_2 R^{2n+2})$ to the Grassmannian at the point $P$, see how to regard it as the 4n-dimensional vector space $\mathrm{Hom}(P, P^\perp)$, and show that a 2n-plane through the origin there is transverse to the bad cone $BC(P)$ if and only if it is the graph of a linear transformation $T: R^{2n} \to R^{2n}$ with no real eigenvalues, with the role of $R^{2n}$ played by two copies of $P^\perp$.

**PROPOSITION 2.** There is a GL(2n, R)-equivariant deformation retraction of the space of linear transformations $T: R^{2n} \to R^{2n}$ with no real eigenvalues to its subspace of linear complex structures $J: R^{2n} \to R^{2n}$.

This is due to Benjamin McKay [2001].

By a **linear complex structure** we mean a linear map $J: R^{2n} \to R^{2n}$ such that $J^2 = -I$. For an **orthogonal complex structure**, we require in addition that the map $J$ be orthogonal.

**PROPOSITION 3.  There is an O(2n)-equivariant deformation retraction of the space of linear complex structures on $R^{2n}$ to its subspace of orthogonal complex structures.**

These results then help us to prove

**PROPOSITION 4.  There exists a smooth fibration $F$ of $S^{2n+1}$ by oriented great circles whose base space $M_F$ is tangent at $P$ to any preassigned 2n-plane transverse to the bad cone $BC(P)$ .**

We then assemble these results to prove theorems A and B .

## CONTENTS

# BACKGROUND

## Nineteenth century: surfaces with simple geodesic behavior.

*Are round two-spheres and real projective planes the only surfaces where all the geodesics are simple closed curves of the same length?* It's easy to understand the motivation for this question: round spheres and real projective planes are the only closed surfaces of constant positive curvature, so it is natural to ask if they are also the only surfaces with such "constant" geodesic behavior. But the answer is *No*, there are other surfaces where all the geodesics are simple closed curves of the same length.

In 1892, Jules Tannery constructed a non-smooth pear-shaped surface in 3-space on which all the geodesics are closed of the same length, except that the equator has half that length. Two years later, Jean Gaston Darboux derived an explicit equation which a surface of revolution in 3-space must satisfy so that all its geodesics are closed, but did not establish the global existence of such a surface. In 1903, Hilbert's student Otto Zoll gave the first example of a smooth surface (in fact, a real analytic surface of revolution) on which all the geodesics are simple closed curves of the same length, other than round spheres and projective planes.

## Twentieth century: Hopf fibrations.

In 1931, Heinz Hopf gave a remarkable example of a map $f$ from the unit 3-sphere $S^3$ to the unit 2-sphere $S^2$. In coordinates:

$$y_1 = 2(x_1x_3 + x_2x_4)$$

$$y_2 = 2(x_2x_3 - x_1x_4)$$

$$y_3 = x_1^2 + x_2^2 - x_3^2 - x_4^2$$

See Figure 3.

This was the first example of a homotopically nontrivial map from a sphere to another sphere of lower dimension, signaling the birth of homotopy theory. Although Hopf presented this map via the above formulas early in his paper, he commented later in that same paper that the great circle fibres of his map were the intersections of the 3-sphere with the complex lines in $C^2$.

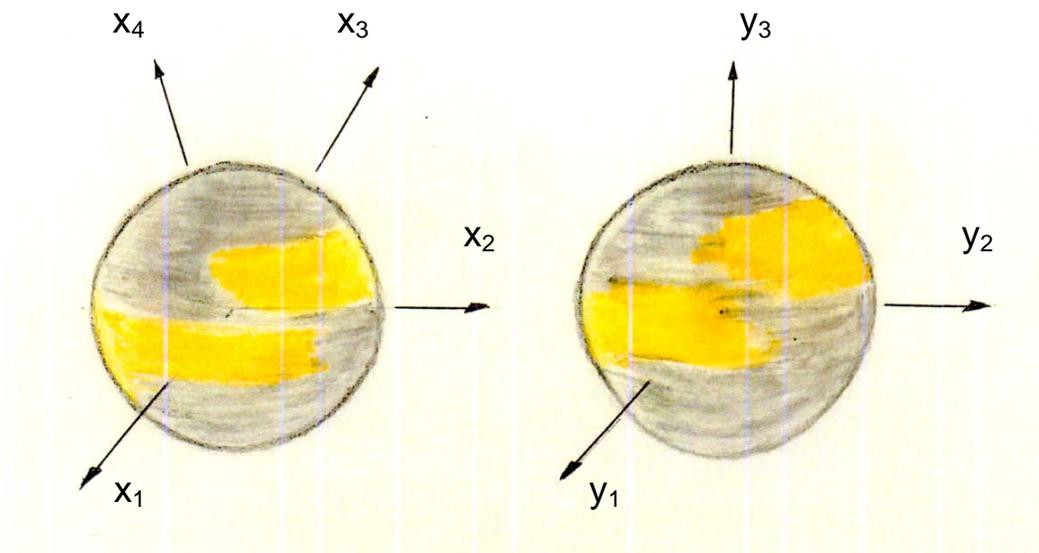

$$S^3 \xrightarrow{\hspace{1.5cm} f \hspace{1.5cm}} S^2$$

**Figure 3. Hopf's map from $S^3$ to $S^2$**

In a follow-up paper in 1935, Hopf presented the higher-dimensional analogues of his first map, using complex numbers, quaternions and Cayley numbers, with the nonassociativity of the Cayley numbers responsible for the truncation of the third series.

$$S^1 \subset S^3 \to S^2 = CP^1, \quad S^1 \subset S^5 \to CP^2, ..., S^1 \subset S^{2n+1} \to CP^n, ...$$

$$S^3 \subset S^7 \to S^4 = HP^1, \quad S^3 \subset S^{11} \to HP^2, ..., S^3 \subset S^{4n+3} \to HP^n, ...$$

$$S^7 \subset S^{15} \to S^8 .$$

## Twentieth century: Blaschke manifolds.

Let M be a closed (compact, no boundary) Riemannian manifold. On each geodesic α from the point p on M, the ***cut point*** is the last point to which α minimizes distance, and the ***cut locus*** C(p) is the set of these.

For example, on a round sphere, the cut locus of each point is just its antipodal point.

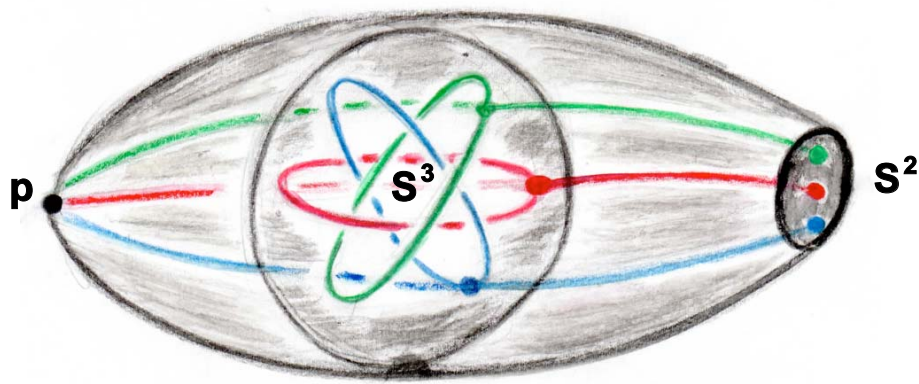

**Figure 4. The complex projective plane CP²**

In this picture of $CP^2$, focus on the point p at the left, and on the geodesics which begin there and eventually coalesce along its cut locus C(p), a round 2-sphere at the right. If we go out along these geodesics any fixed intermediate distance, we come to a 3-sphere on which we record that a circle's worth of geodesics from p will coalesce along each point of C(p). If this intermediate distance is very small, then the 3-sphere is almost round, and its fibration by these circles is almost a Hopf fibration. But as the 3-sphere moves towards the cut locus at the right, these circles will eventually shrink until in the limit they become points, and the 3-sphere collapses to a 2-sphere. The complex projective plane itself is homeomorphic to the mapping cone of this collapsing map $S^3 \to S^2$.

Given the closed Riemannian manifold $M$, if the distance from $p$ to its cut point along $\alpha$ depends neither on the choice of $\alpha$ nor on the choice of $p$, then $M$ is called a ***Blaschke manifold***, the term coined by Marcel Berger [1978].

Examples of Blaschke manifolds are the standard spheres and projective spaces $S^n$, $RP^n$, $CP^n$, $HP^n$ and $CaP^2$, on which all the geodesics from any point come together again after the same distance, independent of direction and point of origin.

The terminology honors Wilhelm Blaschke, who asked, in the first edition [1921] of his *Vorlesungen über Differentialgeometrie*, whether such a surface must be isometric to a round $S^2$ or round $RP^2$.

Reidemeister thought he had a positive answer to Blaschke's question, and this appeared in an appendix to the second edition [1924] of Blaschke's text, but Blaschke pointed out the error in his third edition [1930]. Finally, in 1963, Leon Green proved that a Blaschke surface can only be a round $S^2$ or $RP^2$.

By 1980, the combined work of Marcel Berger, Jerry Kazdan, Alan Weinstein and C.T. Yang showed that Blaschke manifolds "modelled on" $S^n$ and $RP^n$ must, up to scale, be isometric to them. Quite a lot is known about Blaschke manifolds in general, but isometry is known in no other cases.

**What is known about Blaschke manifolds?**

Every Blaschke manifold has the same cohomology ring as one of the spheres or projective spaces mentioned earlier,

$$S^n, \ RP^n, \ CP^n, \ HP^n \ \text{and} \ CaP^2,$$

thanks to the work of Bott [1954] and Samelson [1963], and we say that the Blaschke manifold is ***modelled*** on that standard space.

Here is a summary, due to Benjamin McKay [2013], of what is known to date.

If a Blaschke manifold is modelled on $S^n$ or $RP^n$, then it is (up to scale) ***isometric*** to that model space, thanks to the work of Berger [1978], Kazdan [1978], Weinstein [1974] and Yang [1980].

If a Blaschke manifold is modelled on $CP^n$, then it is ***diffeomorphic*** to this model space, thanks to the work of Yang [1990, 1991] and McKay [2001].

If a Blaschke manifold is modelled on $HP^2$, then it is ***PL-homeomorphic*** to this model space, thanks to the work of Sato and Mizutani [1984].

If a Blaschke manifold is modelled on $HP^n$, then it is ***homotopy equivalent*** to this model space, thanks to the work of Sato [1986].

If a Blaschke manifold is modelled on $CaP^2$, then it is ***homeomorphic*** to this model space, thanks to the work of Gluck-Warner-Yang [1983].

## How do Blaschke manifolds determine fibrations of round spheres by great subspheres?

Let $M$ be a Blaschke manifold, $p$ a point of $M$, and $T_pM$ the tangent space to $M$ at $p$. Let $B(p)$ denote a round ball of radius $r$ in $T_pM$, where $r$ is the common distance from each point of $M$ to its cut locus. See Figure 5.

**Theorem (Omori 1968 and Nakagawa-Shiohama 1970). If $M$ is a Blaschke manifold, then the cut locus $C(p)$ to any point $p$ in $M$ is a smooth submanifold of $M$, and $\exp_p : \partial B(p) \to C(p)$ is a smooth fibre bundle with great subsphere fibres.**

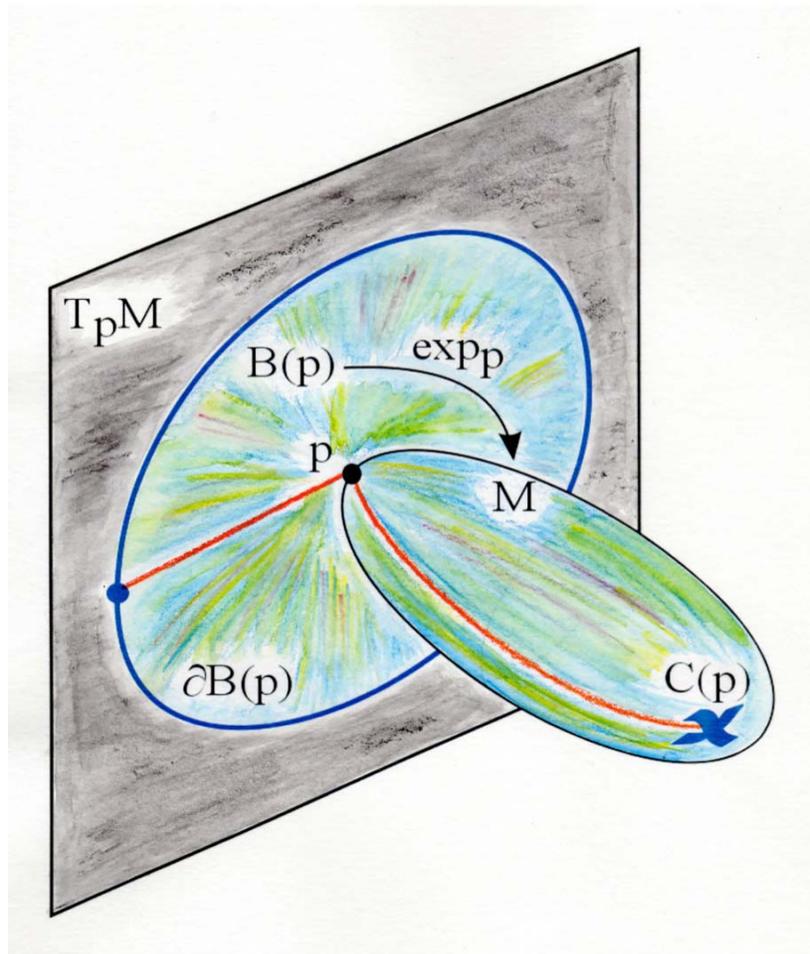

**Figure 5. The exponential map $\exp_p : T_pM \to M$ takes a round ball $B(p)$ onto the Blaschke manifold $M$ and takes $\partial B(p)$ to the cut locus $C(p)$.**

By the above theorem, any Blaschke manifold leads to a smooth fibration of a round sphere by great subspheres. The Blaschke manifold $M$ can be recovered topologically from the fibration $\exp_p : \partial B(p) \to C(p)$, since $M$ is homeomorphic to its mapping cone. Thus to understand Blaschke manifolds topologically, one should understand the topological classification of fibrations of spheres by great subspheres.

**Conjecture.** Any smooth fibration of a sphere by great subspheres is topologically equivalent to a Hopf fibration.

**Caution.** There are many inequivalent fibrations of $S^7$ by 3-spheres (Milnor 1956, Eells and Kuiper 1962), but in general their fibres are not great 3-spheres.

To prove the conjecture, one must figure out how to capitalize on the hypothesis of ***great*** sphere fibres.

The conjecture is known in the following cases:

• Any fibration of $S^3$ by simple closed curves is topologically equivalent to the Hopf fibration [Steenrod 1951].

• Any smooth fibration of $S^7$ by great 3-spheres or of $S^{15}$ by great 7-spheres is topologically equivalent to a Hopf fibration [Gluck-Warner-Yang 1983].

• Any smooth fibration of $S^{2n+1}$ by great circles is smoothly equivalent to a Hopf fibration [Yang 1990 and McKay 2001].

## We know a lot about fibrations of the three-sphere by great circles.

In an old paper with Frank Warner [1983], we studied the different ways in which the three-sphere can be fibered by great circles. We started with a fibration $F$ of $S^3$ by oriented great circles, viewed the base space $M_F$ as a submanifold of the Grassmannian $G_2R^4$ of oriented 2-planes through the origin in $R^4$, used the fact that $G_2R^4$ is isometric to the product of a pair of round 2-spheres, and obtained the following results.

**Theorem I.** A submanifold of $G_2R^4 \cong S^2 \times S^2$ is the base space of a fibration $F$ of $S^3$ by oriented great circles if and only if it is the graph of a distance-decreasing map $f$ from either $S^2$ factor to the other.

**Theorem II.** The great circle fibration $F$ is differentiable if and only if the corresponding distance-decreasing map $f$ is differentiable with $|df| < 1$.

**Theorem III.** Any fibration of $S^3$ by great circles must contain some orthogonal pair of circles.

**Theorem IV.** The space of all oriented great circle fibrations of $S^3$ deformation retracts to the subspace of Hopf fibrations, and hence has the homotopy type of a pair of disjoint two-spheres.

The proofs of these theorems depended crucially on the well known fact that the Grassmannian $G_2R^4$ is isometric to the product of a pair of round two-spheres, and on our introduction of a moduli space for the family of fibrations of $S^3$ by oriented great circles, namely two copies of the set of distance-decreasing maps from $S^2$ to $S^2$.

A correspondingly clear view of the higher Grassmannians is sadly missing from the literature, and even $G_2R^6$, the next one of interest for studying great circle fibrations of $S^5$, seems to be not yet well enough understood to help us find a moduli space for these fibrations.

**But we wish we knew more.** Since we know that any smooth fibration of $S^{2n+1}$ by great circles is smoothly equivalent to a Hopf fibration, we might hope to prove that the set of all such fibrations of $S^{2n+1}$ deformation retracts to its subset of Hopf fibrations. But at present we only know this for $S^3$. Hence the current paper, in which we prove an infinitesimal version of this theorem for great circle fibrations of $S^{2n+1}$.

# THE GRASSMANN MANIFOLD

## Coordinates in the Grassmann manifold $G_2 R^{2n+2}$ .

Given a fibration $F$ of $S^{2n+1}$ by oriented great circles, each fibre $P$ of $F$ lies in and orients some 2-plane through the origin in $R^{2n+2}$ , which we denote by $P$ as well, and so appears as a single point in the Grassmann manifold $G_2 R^{2n+2}$ of all such oriented 2-planes.

The base space $M_F$ of $F$ then appears as a 2n-dimensional topological submanifold of $G_2 R^{2n+2}$ , and if the fibration $F$ is smooth, then the submanifold $M_F$ is also smooth.

Let $P$ be an oriented great circle on $S^{2n+1}$ , equivalently, an oriented 2-plane through the origin in $R^{2n+2}$ , and let $P^\perp$ be its orthogonal complement.

The 4n-dimensional vector space $Hom(P, P^\perp)$ will serve simultaneously as a large coordinate neighborhood about $P$ in $G_2 R^{2n+2}$ , and as the tangent space $T_P(G_2 R^{2n+2})$ to this Grassmann manifold at $P$ , as follows.

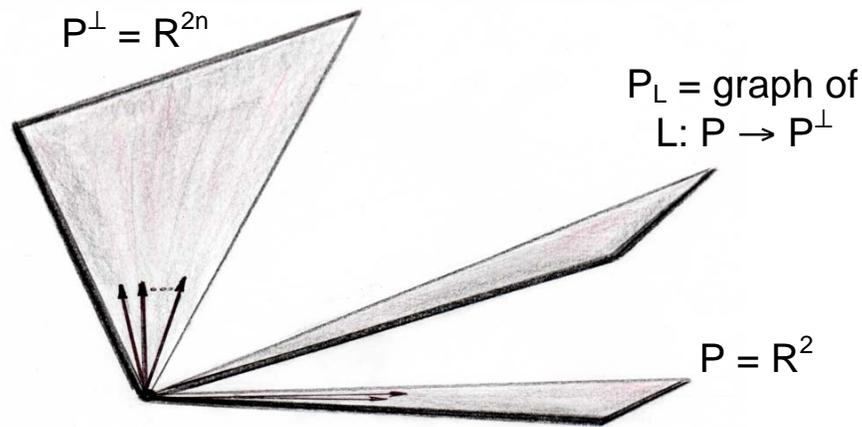

**Figure 6.  $P_L$  =  graph of  $L: P \to P^\perp$  in  $P + P^\perp = R^{2n+2}$**

Suppose that the oriented 2-plane $Q$ in $R^{2n+2}$ contains no vector orthogonal to $P$ , and suppose that its orthogonal projection to $P$ is orientation-preserving. Let $N(P)$ be the collection of all such 2-planes $Q$ . This set $N(P)$ is the domain of our coordinate chart

$$G_2 R^{2n+2} \supset N(P) \longrightarrow \phi \longrightarrow Hom(P, P^\perp) ,$$

defined as follows.

Given $Q \in N(P)$, we can view $Q$ as the graph of a linear transformation $L_Q : P \to P^\perp$ and we set $\phi(Q) = L_Q$. Note that $P$ is itself the graph of the zero transformation, so $\phi(P) = 0$.

Conversely, given a linear transformation $L : P \to P^\perp$, the graph of $L$ is a 2-plane $P_L$ in $R^{2n+2}$, which we may orient via orthogonal projection back to $P$, allowing us to view this graph as an element of $N(P)$.

Since $\text{Hom}(P, P^\perp)$ is a vector space, the differential $\phi_*$ of $\phi : N(P) \to \text{Hom}(P, P^\perp)$ is an isomorphism of the tangent space $T_P(G_2R^{2n+2})$ with $\text{Hom}(P, P^\perp)$.

Thus we may view $\text{Hom}(P, P^\perp)$ simultaneously as a coordinate neighborhood of $P$ in $G_2R^{2n+2}$, and as the tangent space $T_P(G_2R^{2n+2})$ to this Grassmannian at $P$. To connect these two roles, we consider the "identity map"

$$I : T_P(G_2R^{2n+2}) = \text{Hom}(P, P^\perp) \to \text{Hom}(P, P^\perp) = N(P) \subset G_2R^{2n+2}.$$

***Caution.*** $I$ is ***not*** the exponential map: it takes some lines through the origin in $T_P(G_2R^{2n+2}) = \text{Hom}(P, P^\perp)$ to geodesics through $P$ in $G_2R^{2n+2}$, with distortion of parametrization, and takes other lines through the origin to non-geodesics through $P$.

Next we fix bases of $P$ and $P^\perp$ in order to write elements of $\text{Hom}(P, P^\perp)$ as $2n \times 2$ matrices. Let $\{ e_1, e_2 \}$ be an orthonormal basis for $P$, consistent with its orientation. Now orient $P^\perp$ so that the orientations on $P$ and $P^\perp$ together give the orientation on $R^{2n+2}$. Finally, choose an orthonormal basis $\{ f_1, f_2, ..., f_{2n} \}$ for $P^\perp$ consistent with its orientation.

We write elements of $\text{Hom}(P, P^\perp)$ as $2n \times 2$ matrices $A = A_1 \mid A_2$, where $A_1$ and $A_2$ are column $2n$-vectors. We see that $\text{Hom}(P, P^\perp)$ is the sum of two copies of $P^\perp$, since we may write

$$\text{Hom}(P, P^\perp) = \{ A_1 \mid A_2 \} = \{A_1 \mid 0\} + \{0 \mid A_2\} = P^\perp + P^\perp,$$

with the identifications

$$P^\perp = \{A_1 \mid 0\} = \{0 \mid A_2\}.$$

Geometrically, the columns $A_1$ and $A_2$ have the following meaning.

Let $P(t)$ be the oriented 2-plane in $R^{2n+2} = P + P^\perp$ spanned by the frame

$$\{ e_1 + t A_1 , \; e_2 + t A_2 \} .$$

For $-\infty < t < \infty$, this gives us a path $t \to P(t)$ in $G_2R^{2n+2}$ which runs within the domain $N(P)$ of our coordinate chart $\phi : N(P) \to \mathrm{Hom}(P, P^\perp)$. The corresponding path in $\mathrm{Hom}(P, P^\perp)$ is the line $t \to t A_1 \mid t A_2$, and the tangent vector to this path at $t = 0$ is

$$A_1 \mid A_2 \in \mathrm{Hom}(P, P^\perp) = T_P(G_2R^{2n+2}) .$$

## The "bad set" and the "bad cone".

Consider oriented great circle fibrations $F$ of $S^{2n+1}$ which contain a fixed great circle fibre $P$. Because the fibres of $F$ are disjoint, the base space $M_F$ certainly cannot also pass through $Q$ in $G_2 R^{2n+2}$ if the corresponding great circles $P$ and $Q$ intersect on $S^{2n+1}$.

This motivates the following definitions.

The **bad set** $BS(P) \subset G_2 R^{2n+2}$ consists of all oriented 2-planes through the origin in $R^{2n+2}$ which meet $P$ in at least a line. If $M_F$ contains the great circle fibre $P$, then $M_F$ intersects the bad set $BS(P)$ only at $P$ and nowhere else.

The **bad cone** $BC(P) \subset T_P(G_2 R^{2n+2})$ is the tangent cone to the bad set at $P$.

Within the coordinate neighborhood $N(P) = \text{Hom}(P, P^\perp)$ of $P$ in $G_2 R^{2n+2}$, the bad set $BS(P)$ consists of linear transformations $L : P \rightarrow P^\perp$ with nontrivial kernel, because the graphs of such linear transformations intersect $P$ in at least a line. Equivalently, these are the $2n \times 2$ matrices $A = A_1 | A_2$ of rank $0$ or $1$. They all have the form

$$A = A_1 \cos t | A_1 \sin t ,$$

where $A_1$ is a column 2n-vector.

We note that, in the $\text{Hom}(P, P^\perp)$ coordinates on $N(P)$, the portion of the bad set within that neighborhood is a union of lines through the origin $0 = \phi(P)$, namely

$$s A = s A_1 \cos t | s A_1 \sin t , \quad \text{with } -\infty < s < \infty .$$

It follows from this that the tangent cone to the bad set at $P$ coincides with this portion of the bad set, that is,

$$I(BC(P)) = BS(P) \cap N(P) .$$

With abuse of language, we may simply write $BC(P) \subset BS(P)$, and view the bad cone at $P$ as a portion of the bad set at $P$.

**Properties of the bad cone.**

**(1)** In the $\text{Hom}(P, P^\perp)$ coordinates on $N(P)$, the bad cone at $Q$ contains the translate of the bad cone at $P$, namely

$$BC(P) + L_Q \subset BC(Q) ,$$

where $L_Q = \phi(Q)$ in our chart $\phi : N(P) \to \text{Hom}(P, P^\perp)$ centered at $P$.

That's because the linear transformations $L_{Q*} : P \to P^\perp$ which correspond to points of $BC(Q)$ are those which agree with $L_Q$ on some nonzero vector $u$ in $P$. Thus $L_{Q*} - L_Q$ contains $u$ in its kernel, and hence belongs to $BC(P)$.

**(2)** The bad cone $BC(P)$ is homeomorphic to a cone over $S^1 \times S^{2n-1}$.

We see this as follows. If $L : P \to P^\perp$ is a linear transformation with a nontrivial kernel, then its $2n \times 2$ matrix $A$ has the form

$$A = \cos t \; A_1 \mid \sin t \; A_1 ,$$

where $A_1$ is some column $2n$-vector.

If we fix $t$ and let $A_1$ vary, we get a $2n$-plane which is part of the bad cone.

If we then let $t$ vary, we fill out the bad cone with a circle's worth of such $2n$-planes, modulo the involution $(t, A_1) \to (t + \pi, -A_1)$.

Equivalently, $BC(P)$ is a cone over the quotient of $S^1 \times S^{2n-1}$ by this involution. But this quotient is homeomorphic to $S^1 \times S^{2n-1}$, since the antipodal map on an odd-dimensional sphere is isotopic to the identity.

In similar fashion, the bad set $BS(P)$ is homeomorphic to the suspension of $S^1 \times S^{2n-1}$.

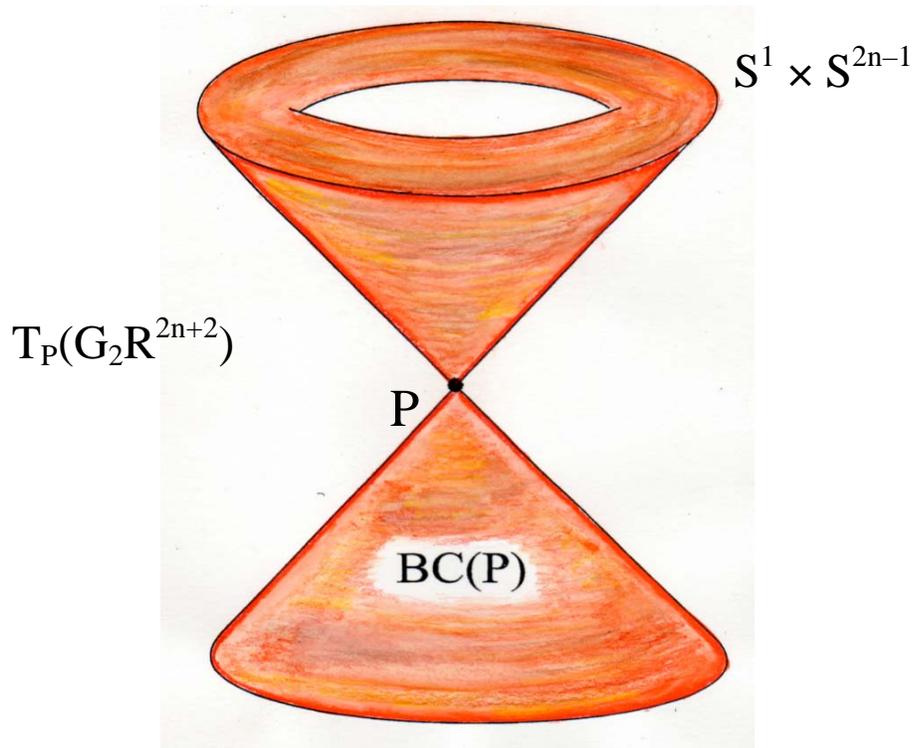

**Figure 7. The bad cone BC(P)**

When we come to Proposition 1, we will visualize the Grassmann manifold $G_2 R^{2n+2}$ with a bad cone $BC(P)$ inside the tangent space $T_P(G_2 R^{2n+2})$ at **each** of its points $P$, thus giving us a **field of bad cones**, as shown in Figure 8.

# PROOF OF PROPOSITION 1

Now we characterize the smooth submanifolds of $G_2 R^{2n+2}$ which correspond to the base space of some smooth fibration of $S^{2n+1}$ by great circles.

**PROPOSITION 1. A closed smooth 2n-dimensional submanifold M of $G_2 R^{2n+2}$ is the base space of a smooth fibration of $S^{2n+1}$ by great circles if and only if it is transverse to the bad cone at each of its points.**

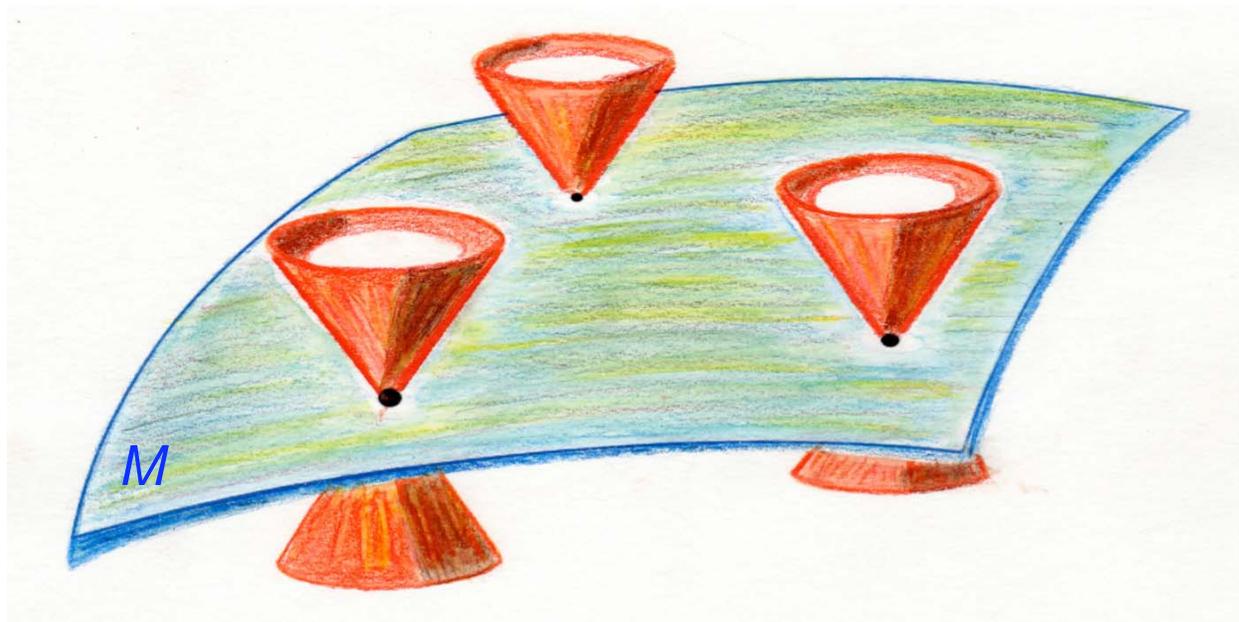

**Figure 8. M is like a submarine negotiating a mine field**

**Proof.**

Suppose first that $F$ is a smooth fibration of $S^{2n+1}$ by oriented great circles. We want to show that its base space $M_F$ in $G_2 R^{2n+2}$ is a smooth submanifold transverse to the field of bad cones there.

For $S^3$, this is Theorem B of [G-W, 1983].

For smooth fibrations of spheres by great subspheres of any dimension, this is Theorem 4.1 of [G-W-Y, 1983].

This was proved again for all great circle fibrations of $S^{2n+1}$ by Benjamin McKay [2001], from a different point of view.

Suppose, conversely, that $M$ is a closed, smooth 2n-dimensional submanifold of $G_2R^{2n+2}$ which is transverse to the field of bad cones.

There is a canonical $S^1$ bundle $E$ over $G_2R^{2n+2}$ whose fibre over $P$ is the great circle in the 2-plane $P$. Let $E_M$ be the restriction of this bundle to the submanifold $M$,

$$E_M = \{(P, v) : P \in M , v \in P , \|v\| = 1 \} .$$

Let $\rho: E_M \to M$ be the projection map, and let $g: E_M \to S^{2n+1}$ be the map which includes each great circle fibre into $S^{2n+1}$, that is, $g(P, v) = v$.

Our task is to show that $g$ is a diffeomorphism.

First, we claim that transversality of $M$ to the bad cone through each of its points implies that the map $g$ is an immersion.

Suppose, to the contrary, that $dg$ has a nontrivial kernel at some point $v$ in the fibre $P$.

Consider a path $\gamma: (-1, 1) \to E_M$, and write $\gamma(t) = (P(t), v(t))$, such that $\gamma(0) = (P, v)$, with $\gamma'(0) \neq 0$.

We will show that if $\gamma'(0)$ is in the kernel of the derivative $dg_v$, then $M$ must be tangent to the bad cone $BC(P)$ at $P$.

Consider the path $P(t) = \rho \gamma(t)$ in $M$, with $P(0) = P$.

Using the coordinate neighborhood $\text{Hom}(P, P^\perp)$ about $P$ in $G_2R^{2n+2}$, the path $P(t)$ corresponds to a path $L(t)$ in $\text{Hom}(P, P^\perp)$.

Since $P(0) = P$, we have $L(0) = 0$.

Now $g\gamma(t) = v(t)$ lies in $P(t)$, which is the graph of $L(t)$, so we can write

$$g\gamma(t) = (w(t) , L(t) w(t))$$

as an ordered pair of vectors in $P \times P^\perp$, with $w(0) \neq 0$.

We differentiate with respect to $t$ and set $t = 0$ to get

$$(g\,\gamma)'(0) = (\,w'(0)\,,\ L'(0)\,w(0) + L(0)\,w'(0)\,) \in P \times P^{\perp}\,.$$

Now we are assuming that $(g\,\gamma)'(0) = 0$ in $R^{2n+2} = P + P^{\perp}$ and we know that $L(0) = 0$, so we conclude that $L'(0)\,w(0) = 0$.

Since $w(0) \neq 0$, this tells us that $L'(0)$ has a nontrivial kernel, and hence lies in the bad cone $BC(P)$ at $P$.

Therefore the path $P(t)$ in $M$ is tangent to the bad cone at $P(0) = P$, contrary to the assumption that $M$ is transverse to the field of bad cones.

So we have just shown that the map $g\colon E_M \to S^{2n+1}$ is an immersion.

But $E_M$ is compact, and so the map $g$ is both open and closed, and hence its image $g(E_M)$ must be all of $S^{2n+1}$.

Thus $g$ is a covering map, and since $S^{2n+1}$ is simply connected for $n \geq 1$, $g$ must be a diffeomorphism.

Thus $E_M$ gives a smooth fibration of $S^{2n+1}$ by great circles, with $M$ as its base space, completing the proof of the lemma.

**Remarks.**

**(1)** The proofs in [G-W, 1983] and in [G-W-Y, 1983] that the base space $M_F$ of a smooth fibration $F$ by great subspheres is transverse to the field of bad cones use the fact that the local trivializations of $F$ are diffeomorphisms.

One can have a topological fibration $F$ of $S^{2n+1}$ by great circles whose base space $M_F$ is a smooth submanifold of $G_2R^{2n+2}$ occasionally tangent to a bad cone, and then the local trivializations of $F$ will be smooth homeomorphisms, but not diffeomorphisms.

**(2)** A small, smooth 2n-disk in $G_2R^{2n+2}$ which is transverse to the field of bad cones gives a fibration of an open tube in $S^{2n+1}$ by great circles.

# PROOF OF PROPOSITION 2

## 2n × 2n matrices with no real eigenvalues.

In this section, we will see how $2n \times 2n$ matrices with no real eigenvalues arise in our study of 2n-planes tangent to the base space of a smooth fibration of $S^{2n+1}$ by great circles.

In the 4n-dimensional vector space $\mathrm{Hom}(P, P^{\perp}) = P^{\perp} + P^{\perp}$, we need to recognize those 2n-dimensional subspaces which are transverse to the bad cone $BC(P)$, since they will be precisely those, according to Propositions 1 and 4, which can serve as tangent spaces to the base spaces of fibrations of $S^{2n+1}$ by great circles.

**LEMMA.** A 2n-dimensional subspace of $\mathrm{Hom}(P, P^{\perp}) = P^{\perp} + P^{\perp}$ is transverse to the bad cone $BC(P)$ if and only if it is the graph of a linear map with no real eigenvalues from one $P^{\perp}$ summand to the other.

**Proof.** A 2n-dimensional subspace $T$ of $\mathrm{Hom}(P, P^{\perp})$ transverse to the bad cone can meet each of the two summands $P^{\perp} = \{A_1 \mid 0\}$ and $P^{\perp} = \{0 \mid A_2\}$ only at the origin, since these summands lie entirely in the bad cone. Hence $T$ is the graph of a linear map $L_T : P^{\perp} \rightarrow P^{\perp}$ between these subspaces, in either order.

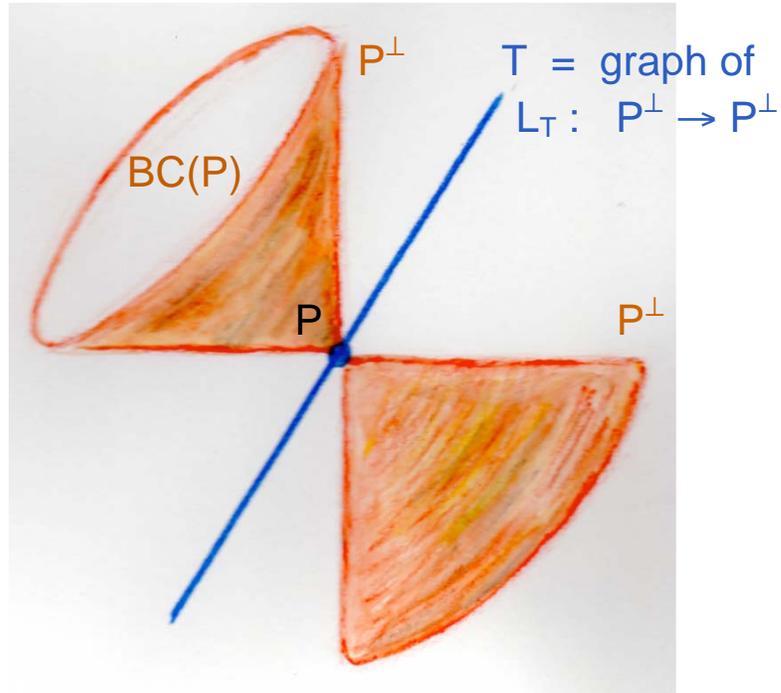

**Figure 9. T is transverse to the bad cone BC(P) if and only if
it is the graph of a linear map $L_T : P^\perp \to P^\perp$ with no real eigenvalues**

If $L_T$ has a real eigenvalue $\lambda$ with eigenvector $A_1$, then its graph $T$ contains
the vector $A_1 \mid \lambda A_1$, a $2n \times 2$ matrix of rank 1, hence in the bad cone $BC(P)$.

Thus a 2n-dimensional subspace $T$ of $\mathrm{Hom}(P, P^\perp)$ which is transverse to the
bad cone is the graph of a linear map $L_T$ as above with no real eigenvalues.

Conversely, if $T$ is a 2n-dimensional subspace of $\mathrm{Hom}(P, P^\perp)$ which is the
graph of a linear map $L_T : P^\perp \to P^\perp$ with no real eigenvalues, then $T$ contains
no $2n \times 2$ matrices of rank 1, and so is transverse to the bad cone $BC(P)$,
proving the Lemma.

## Improving maps with no real eigenvalues.

Recall that by a ***linear complex structure*** we mean a linear map $J: R^{2n+2} \to R^{2n+2}$ such that $J^2 = -I$, and that for an ***orthogonal complex structure***, we require in addition that the map $J$ be orthogonal.

Given any orthogonal complex structure $J: R^{2n+2} \to R^{2n+2}$, the unit circles in the J-complex lines yield a Hopf fibration $H$ of $S^{2n+1}$ by oriented great circles.

**LEMMA.** The tangent 2n-plane to the base space $M_H$ at a complex line $P$ is the graph of $J|_{P^\perp} : P^\perp \to P^\perp$.

**Proof.**

The points $L$ in the large coordinate neighborhood $Hom(P, P^\perp)$ of $P$ in $G_2 R^{2n+2}$ are represented by $2n \times 2$ matrices $A = A_1 | A_2$, where the two columns are the L-images in $P^\perp$ of an ON basis $e_1, e_2$ for $P$ with $J(e_1) = e_2$.

The points $Q$ in this neighborhood which lie in the base space $M_H$ of the fibration $H$ are J-complex lines, meaning images of a J-complex linear map $L: P \to P^\perp$. Since $L(e_1) = A_1$ and $L(e_2) = A_2$ and $L \circ J = J \circ L$, we have

$$A_2 = L(e_2) = L(J(e_1)) = J(L(e_1)) = J(A_1).$$

Thus the points of $M_H$ in this coordinate neighborhood lie on the graph of $J|_{P^\perp} : P^\perp \to P^\perp$.

Since the coordinate neighborhood $Hom(P, P^\perp)$ of $P$ serves as its own tangent space at $P$, the graph of $J|_{P^\perp} : P^\perp \to P^\perp$ serves as the tangent 2n-plane to $M_H$ at $P$, as claimed.

**Remarks.**

(1) We note that the portion of $M_H$ within the large open neighborhood $Hom(P, P^\perp)$ of $P$ in $G_2 R^{2n+2}$ appears as a 2n-plane through the origin there.

(2) The above Lemma and Remark hold equally well if $J: R^{2n+2} \to R^{2n+2}$ is only a linear complex structure such that $J(P^\perp) = P^\perp$.

**PROPOSITION 2.  There is a GL(2n, R)-equivariant deformation retraction of the space of linear transformations  T: $R^{2n} \to R^{2n}$  with no real eigenvalues to its subspace of linear complex structures J: $R^{2n} \to R^{2n}$ .**

**Proof.**

See Benjamin McKay [2001], pages 16 - 20 .

Let  T: $R^{2n} \to R^{2n}$  be a linear transformation with no real eigenvalues. Complexify  $R^{2n}$  to get  $C^{2n}$ , and regard  T: $C^{2n} \to C^{2n}$ . Since  T  is real, its eigenvalues  $\lambda$  occur in conjugate pairs.

Split  $C^{2n}$  into a direct sum  $\sum_\lambda E_\lambda T$  of the generalized eigenspaces of  T , where

$$E_\lambda T \;=\; \{\, v \in C^{2n} : (T - \lambda I)^k \, v \;=\; 0 \;\text{ for some }\; k > 0 \,\} \;,$$

with dim($E_\lambda T$) = multiplicity of the eigenvalue  $\lambda$ .  Complex conjugation in  $C^{2n}$  takes  $E_\lambda T$  to  $E_{\bar\lambda} T$  since  T  is real.

Reorganize the direct sum,

$$C^{2n} \;=\; \sum_{\text{Im}\,\lambda > 0} E_\lambda T \;+\; \sum_{\text{Im}\,\lambda < 0} E_\lambda T \;=\; V_C^{\,+} \;+\; V_C^{\,-} \;,$$

and note that complex conjugation interchanges  $V_C^{\,+}$  and  $V_C^{\,-}$ .

Now define a complex linear map  $J_T$: $C^{2n} \to C^{2n}$  by  $J_T(v) = i\,v$ if  $v \in V_C^{\,+}$  and  $J_T(v) = -i\,v$  if  $v \in V_C^{\,-}$ .  This map  $J_T$  commutes with complex conjugation, and hence takes real vectors to real vectors, so that  $J_T$: $R^{2n} \to R^{2n}$  is a linear complex structure.

It is clear from construction that the correspondence  $T \to J_T$  is GL(2n, R)-equivariant.

Our desired deformation retraction is given by the formula

$$T_t \;=\; (1-t)\,T \;+\; t\,J_T\,.$$

One easily checks by looking at the blocks in the Jordan normal form for $T$ that each of the transformations $T_t$ has no real eigenvalues.

Since $T$ and $J_T$ each commute with complex conjugation, the same is true for $T_t$, and hence it also takes real vectors to real vectors.

To confirm that the proposed deformation retraction $T_t$ depends continuously on $T$, we must check that $J_T$ itself depends continuously on $T$.

Since $J_T$ is defined as multiplication by $i$ on $V_C^{+}$ and by $-i$ on $V_C^{-}$, this amounts to checking that the subspaces $V_C^{+}$ and $V_C^{-}$ depend continuously on the choice of $T$ from among the linear transformations $R^{2n} \to R^{2n}$ with no real eigenvalues.

This is implied by Lemma 6 on page 18 of Benjamin McKay [2001], where he shows that the map $T \to J_T$ is the projection of a smooth fibre bundle.

We give a different argument here.

Let $\lambda_1$, $\lambda_2$, ..., $\lambda_n$ be the eigenvalues of $T$ with positive imaginary part, and $\bar{\lambda}_1$, $\bar{\lambda}_2$, ..., $\bar{\lambda}_n$ their complex conjugates, which are the eigenvalues of $T$ with negative imaginary part. In each case, an eigenvalue may be listed several times according to its multiplicity.

Consider the complex polynomials

$$p_T^{+}(z) \;=\; (z-\lambda_1)\,(z-\lambda_2)\,...\,(z-\lambda_n) \;\text{ and }\; p_T^{-}(z) \;=\; (z-\bar{\lambda}_1)\,(z-\bar{\lambda}_2)\,...\,(z-\bar{\lambda}_n)\,,$$

which are the characteristic polynomials of the restrictions of $T$ to $V_C^{+}$ and $V_C^{-}$, respectively. Their product $p_T(z) = p_T^{+}(z)\,p_T^{-}(z)$ is the characteristic polynomial of $T$ on all of $V_C$.

By the Cayley-Hamilton theorem, the linear transformation $p_T{}^+(T)$ vanishes on $V_C{}^+$, the linear transformation $p_T{}^-(T)$ vanishes on $V_C{}^-$, while their product (composition) $p_T(T) = p_T{}^+(T) \, p_T{}^-(T)$ vanishes on all of $V_C$.

Since $p_T{}^+(z)$ and $p_T{}^-(z)$ have no roots in common, they are relatively prime, and hence there are polynomials $a_T{}^+(z)$ and $a_T{}^-(z)$ such that

$$a_T{}^+(z) \, p_T{}^+(z) \; + \; a_T{}^-(z) \, p_T{}^-(z) \; = \; 1 \, .$$

Inserting $T$ in place of $z$, we get

(*) $$a_T{}^+(T) \, p_T{}^+(T) \; + \; a_T{}^-(T) \, p_T{}^-(T) \; = \; I \, .$$

**LEMMA.** The kernels of the linear maps $p_T{}^+(T)$ and $p_T{}^-(T) : V_C \to V_C$ are precisely

$$\ker p_T{}^+(T) \; = \; V_C{}^+ \quad \text{and} \quad \ker p_T{}^-(T) \; = \; V_C{}^- \, .$$

**Proof.** We already know that $p_T{}^+(T)$ vanishes on $V_C{}^+$, so that $\ker p_T{}^+(T)$ contains $V_C{}^+$, and likewise $\ker p_T{}^-(T)$ contains $V_C{}^-$. Now $V_C{}^+$ and $V_C{}^-$ are complex n-dimensional subspaces of the complex 2n-dimensional space $V_C$. If either $\ker p_T{}^+(T)$ is larger than $V_C{}^+$ or $\ker p_T{}^-(T)$ is larger than $V_C{}^-$, then there would have to be a nonzero vector $v$ in $V_C$ which lies in both kernels. But then applying formula (*) above to $v$ would give a contradiction, because the left side would kill $v$, while the right side would preserve it. This completes the proof of the lemma.

Now as $T$ varies continuously among linear transformations $R^{2n} \to R^{2n}$ with no real eigenvalues, the roots of its characteristic polynomial also vary continuously (with multiple roots permitted to split into simpler ones), and so by the above lemma, the subspaces $V_C{}^+$ and $V_C{}^-$ also vary continuously.

This completes the proof of Proposition 2.

# PROOF OF PROPOSITION 3

Now we discuss the second step of our deformation retraction.

**PROPOSITION 3.  There is an O(2n)-equivariant deformation retraction of the space of linear complex structures on $R^{2n}$ to its subspace of orthogonal complex structures.**

To prove this, we will use the one-to-one correspondence between linear complex structures $J: R^{2n} \to R^{2n}$ and direct sum decompositions of $C^{2n} = V_C^+ + V_C^-$ into a pair of conjugate complex subspaces, the $+i$ and $-i$ eigenspaces of $J$ on $C^{2n}$, as described in the proof of Proposition 2.

We will check that the complex structure $J$ is orthogonal if and only if $V_C^+$ and $V_C^-$ are orthogonal to one another.

Our goal will then be to describe a deformation retraction from the set of all pairs $V_C^+$ and $V_C^-$ of complex n-dimensional conjugate subspaces of $C^{2n}$ to its subset of orthogonal such pairs. Intuitively, this deformation retraction is given by opening up all the angles between $V_C^+$ and $V_C^-$ in a coordinated fashion until they become orthogonal.

We turn now to providing the details.

**Characterization of orthogonal complex structures.**

**LEMMA.** A linear complex structure $J: R^{2n} \to R^{2n}$ is orthogonal if and only if $v$ and $J(v)$ are orthogonal to one another for all vectors $v$ in $R^{2n}$.

**Proof.**

If $J$ is an orthogonal complex structure, it is easy to check that $v$ and $J(v)$ are orthogonal to one another for all vectors $v$ in $R^{2n}$.

In the other direction, suppose that $J: R^{2n} \to R^{2n}$ is a linear complex structure for which $v$ and $J(v)$ are orthogonal for all vectors $v$ in $R^{2n}$.

Apply this statement to the vector $w = u + J(v)$ to learn that

$$0 = w \bullet J(w) = (u + J(v)) \bullet J(u + J(v))$$

$$= (u + J(v)) \bullet (J(u) + J^2(v)) = (u + J(v)) \bullet (J(u) - v)$$

$$= u \bullet J(u) - u \bullet v + J(v) \bullet J(u) - J(v) \bullet v$$

$$= - u \bullet v + J(u) \bullet J(v),$$

from which we get $J(u) \bullet J(v) = u \bullet v$, confirming that $J$ is an orthogonal transformation.

**LEMMA.** A linear complex structure $J: R^{2n} \to R^{2n}$ is orthogonal if and only if the conjugate complex subspaces $V_C^+$ and $V_C^-$ of $C^{2n}$ are orthogonal to one another.

**Proof.**

We start with $R^2$, and let $J: R^2 \to R^2$ be given by the matrix $\begin{matrix} 0 & b \\ -1/b & 0 \end{matrix}$.

It is easy to see by continuity that every complex structure $J$ on $R^2$ moves some nonzero vector orthogonal to itself, so that it can be expressed in the above matrix form for some orthonormal basis.

The above map $J$ is orthogonal if and only if $b = \pm 1$.

The eigenvalues of $J$ are $i$ and $-i$, and corresponding eigenvectors of $J$ on $C^2$ are the column vectors $u = [b \ \ i]$ and $v = [b \ \ -i]$.

The complex subspaces $V_C^+$ and $V_C^-$ of $C^2$ are generated in this case by the $i$ and $-i$ eigenvectors above. That is,

$\quad V_C^+ = C\{u = [b \ \ i]\} = R\{u = [b \ \ i] , \ u' = iu = [ib \ \ -1]\}$ and

$\quad V_C^- = C\{v = [b \ \ -i]\} = R\{v = [b \ \ -i] , \ v' = iv = [ib \ \ 1]\}$.

We compute the dot products of these vectors and learn that

$\quad\quad u \bullet v = b^2 - 1 , \quad u \bullet v' = 0 , \quad u' \bullet v = 0 , \quad u' \bullet v' = b^2 - 1 .$

Hence the $+i$ and $-i$ eigenspaces $V_C^+$ and $V_C^-$ are orthogonal to one another if and only if $b = \pm 1$, which is precisely the condition that the complex structure $J$ be orthogonal.

This completes the argument for $R^2$.

With this in hand, we carry out the general argument for $R^{2n}$.

If $J : R^{2n} \rightarrow R^{2n}$ is an orthogonal complex structure, then we can choose an orthonormal basis for $R^{2n}$ with respect to which the matrix for $J$ is in block diagonal form, with $2 \times 2$ blocks

$$\begin{array}{cc} 0 & -1 \\ 1 & 0 \end{array}$$

down the diagonal.

Then $V_C^+$ and $V_C^-$ are each complex n-dimensional subspaces of $C^{2n}$. Each is an orthogonal direct sum of complex lines. The $r^{th}$ complex lines in each direct sum are orthogonal to one another by the completed task in $R^2$, whereas the $r^{th}$ complex line in one sum is automatically orthogonal to the $s^{th}$ complex line in the other sum when $r \neq s$. It follows that the complex subspaces $V_C^+$ and $V_C^-$ are othogonal to one another in $C^{2n}$.

If $J : R^{2n} \rightarrow R^{2n}$ is **not** an orthogonal complex structure, then it follows from our earlier Lemma characterizing orthogonal complex structures that there is some vector $v$ in $R^{2n}$ for which $J(v)$ is not orthogonal to $v$.

The 2-plane spanned by this $v$ and $J(v)$ is invariant under $J$, but on it $J$ is **not** a rotation by $90^o$, as we saw in $R^2$, and hence $V_C^+$ and $V_C^-$ are **not** orthogonal to one another.

This completes the proof of the Lemma.

# Principal angles.

We discuss the notion of ***principal angles*** in three settings:

> (1) between a pair of real linear subspaces in $R^n$ ,

> (2) between a pair of complex linear subspaces in $C^n$ ,

> (3) between a complex linear subspace and
>    its complex conjugate subspace in $C^{2n}$ .

The intention is to characterize the ***relative position*** of the two subspaces, up to the action of an appropriate group of isometries of the ambient space, which in the three cases above are the groups $O(n)$ , $U(n)$ , and $O(2n)$ .

The notion and use of principal angles in the real setting (1) is familiar in geometry, and goes back at least to Camille Jordan [1875]; see also Gluck [1967]. But the extension to the complex settings (2) and (3) appears to be much less familiar, though we note the papers by Scharnhorst [2001] and by Galantai and Hegedus [2006], the latter having a very nice set of references.

## (1) Principal angles between a pair of linear subspaces in $R^n$ .

Let $P$ and $Q$ be k-planes through the origin in $R^n$ . Then the relative position of $P$ and $Q$ in $R^n$ is characterized up to the action of $O(n)$ by $k$ principal angles $\theta_1$ , $\theta_2$ , ... , $\theta_k$ , obtained as follows.

$\theta_1$ is the smallest angle that any vector in $P$ makes with any vector in $Q$ . Pick such unit vectors $v_1$ in $P$ and $w_1$ in $Q$ . Let $P_2$ be the orthogonal complement of $v_1$ in $P = P_1$ and let $Q_2$ be the orthogonal complement of $w_1$ in $Q = Q_1$ . Thus $P_2$ and $Q_2$ are k−1 planes through the origin in $R^n$ .

***Remark.*** It follows easily from the minimality of $\theta_1$ that $P_2$ is also orthogonal to $w_1$ , and that $Q_2$ is also orthogonal to $v_1$ .

We move to the induction step. If $\theta_1 = 0$ , then $v_1 = w_1$ and we replace $R^n$ by the $R^{n-1}$ orthogonal to $v_1 = w_1$ , and replace the k-planes $P$ and $Q$ by the k−1 planes $P_2$ and $Q_2$ .

If $\theta_1 > 0$, then $v_1$ and $w_1$ are independent and span a 2-plane through the origin. We replace $R^n$ by the $R^{n-2}$ orthogonal to this 2-plane, and replace the k-planes $P$ and $Q$ by the k−1 planes $P_2$ and $Q_2$. In this case we need the above remark, to guarantee that $P_2$ and $Q_2$ lie in this $R^{n-2}$.

Now we iterate the construction, with $R^n$ replaced by either $R^{n-1}$ or $R^{n-2}$ as detailed above, and with $P$ and $Q$ replaced by $P_2$ and $Q_2$.

Following through to the end, we get orthonormal bases

$$v_1, v_2, \dots, v_k \quad \text{and} \quad w_1, w_2, \dots, w_k$$

for the k-planes $P$ and $Q$, respectively, with ***principal angles***

$$\theta_1 \leq \theta_2 \leq \dots \leq \theta_k \leq \pi/2$$

between the vectors $v_1$ and $w_1$, $v_2$ and $w_2$, $\dots$, $v_k$ and $w_k$, and with $v_r$ orthogonal to $w_s$ for $r \neq s$.

The principal angles between $P$ and $Q$ characterize their relative position in $R^n$ as follows.

**(1) PRINCIPAL ANGLES THEOREM IN $R^n$.** Let $P$ and $Q$ be a pair of k-planes through the origin in $R^n$, and likewise for $P'$ and $Q'$. Then there is a rigid motion (element of O(n)) taking $P$ to $P'$ and simultaneously taking $Q$ to $Q'$ if and only if the principal angles between $P$ and $Q$ are the same as those between $P'$ and $Q'$.

**Proof.**

The condition of matching principal angles is clearly necessary for the existence of such a rigid motion.

Conversely, if the principal angles $\theta_1 \leq \theta_2 \leq \dots \leq \theta_k$ between $P$ and $Q$ match the principal angles $\theta'_1 \leq \theta'_2 \leq \dots \leq \theta'_k$ between $P'$ and $Q'$, then we easily obtain a rigid motion of $R^n$ which takes the orthonormal bases $v_1, v_2, \dots, v_k$ and $w_1, w_2, \dots, w_k$ for $P$ and $Q$ to the orthonormal bases $v'_1, v'_2, \dots, v'_k$ and $w'_1, w'_2, \dots, w'_k$ for $P'$ and $Q'$.

**(2) Principal angles between a pair of complex linear subspaces of $C^n$.**

Let $P$ and $Q$ be complex k-dimensional linear subspaces of $C^n$, which to real eyes look like 2k-planes through the origin in $R^{2n}$.

To get principal angles between $P$ and $Q$, and corresponding orthonormal bases for each of them, we begin as in the real case. Let $\theta_1$ be the smallest angle that any vector in $P$ makes with any vector in $Q$, and pick such unit vectors $v_1$ in $P$ and $w_1$ in $Q$.

Then consider $i\,v_1$ and $i\,w_1$. These will be another pair of unit vectors in $P$ and $Q$, respectively, since each of these is a **complex** linear subspace. The angle between $i\,v_1$ and $i\,w_1$ is also $\theta_1$, because multiplication by $i$ is an isometry of $C^n$ which takes $P$ to itself and $Q$ to itself.

The list of principal angles begins with $\theta_1$, $\theta_1$, while our orthonormal bases for $P$ and $Q$ over the reals begin with $v_1$, $i\,v_1$ for $P$ and $w_1$, $i\,w_1$ for $Q$.

We economize and list angles and bases from a complex point of view, so that our principal angles begin with just $\theta_1$, while our orthonormal bases for $P$ and $Q$ over the complex numbers begins with $v_1$ for $P$ and $w_1$ for $Q$.

We then iterate, as in the real case, and end with complex orthonormal bases

$$v_1,\ v_2,\ \dots,\ v_k \quad \text{and} \quad w_1,\ w_2,\ \dots,\ w_k$$

for the k-planes $P$ and $Q$, with **principal angles**

$$\theta_1 \leq \theta_2 \leq \dots \leq \theta_k \leq \pi/2$$

between the vectors $v_1$ and $w_1$, $v_2$ and $w_2$, $\dots$, $v_k$ and $w_k$, and with $v_r$ orthogonal to $w_s$ for $r \neq s$.

**(2) PRINCIPAL ANGLES THEOREM IN $C^n$.** Let $P$ and $Q$ be a pair of complex k-planes through the origin in $C^n$, and likewise for $P'$ and $Q'$. Then there is an element of $U(n)$ taking $P$ to $P'$ and simultaneously taking $Q$ to $Q'$ if and only if the principal angles between $P$ and $Q$ are the same as those between $P'$ and $Q'$.

We omit the proof, which is basically the same as in the real case.

**(3) Principal angles between conjugate complex linear subspaces in $C^{2n}$ .**

Let $P^k$ and $\overline{P}^k$ be conjugate complex subspaces of $C^{2n}$ which meet only at the origin. We want to define the principal angles between them.

Let $\theta_1$ be the smallest angle that any complex line $L$ in $P^k$ makes with its conjugate complex line $\overline{L}$ in $\overline{P}^k$ . We claim that there will be a unit vector $v_1$ in $L$ which makes that angle $\theta_1$ with its complex conjugate $\overline{v}_1$ in $\overline{L}$ .

The reason for this is that the nearest neighbor map from the unit circle in $L$ to the unit circle in $\overline{L}$ is orientation-preserving, while the complex conjugation map between these unit circles is orientation-reversing. So there is sure to be a coincidence between these two maps, meaning a unit vector $v_1$ in $L$ whose nearest neighbor in $\overline{L}$ is its own conjugate $\overline{v}_1$ .

Thus $v_1$ makes the angle $\theta_1$ with $\overline{v}_1$ , and likewise $i\,v_1$ makes that same angle $\theta_1$ with $i\,\overline{v}_1$ . We note that $i\,v_1$ and $i\,\overline{v}_1$ , though nearest neighbors in $L$ and $\overline{L}$ , are **not** complex conjugates of one another.

Now let $P_2$ be the orthogonal complement of the complex line $L = C\,v_1$ in $P^k$ , and then $\overline{P}_2$ will automatically be the orthogonal complement of the complex line $\overline{L} = C\,\overline{v}_1$ in $\overline{P}^k$ .

***Remark.*** As in the previous two cases, we find that $P_2$ is also orthogonal to $\overline{L} = C\,\overline{v}_1$ , and then (automatically) $\overline{P}_2$ is also orthogonal to $L = C\,v_1$ , and omit the details.

Then, since $P^k$ and $\overline{P}^k$ meet only at the origin, we have $\theta_1 > 0$ .

So we replace $C^{2n}$ by the $C^{2n-2}$ orthogonal to $C\,v_1 + C\,\overline{v}_1$ , and replace $P$ and $\overline{P}$ by the complex $k-1$ dimensional subspaces $P_2$ and $\overline{P}_2$ , both lying in this $C^{2n-2}$ , thanks to the above remark.

As before, we iterate the construction, with $C^{2n}$ replaced by $C^{2n-2}$ and with $P$ and $\overline{P}$ replaced by $P_2$ and $\overline{P}_2$ .

Following through to the end, we get complex orthonormal bases

$$v_1, \ v_2, \ ..., v_k \quad \text{and} \quad \overline{v}_1, \ \overline{v}_2, \ ..., \ \overline{v}_k$$

for the k-planes $P^k$ and $\overline{P}^k$, respectively, with ***constrained principal angles***

$$0 < \theta_1 \leq \theta_2 \leq ... \leq \theta_k \leq \pi/2$$

between the vectors $v_1$ and $\overline{v}_1$, $v_2$ and $\overline{v}_2$, ..., $v_k$ and $\overline{v}_k$, and with $C\,v_r$ orthogonal to $C\,\overline{v}_s$ for $r \neq s$.

***Remark.*** The "constraint" on these principal angles is seen at the beginning, when we minimize the angle $\theta_1$ between a complex line $L$ in $P^k$ and its conjugate $\overline{L}$ in $\overline{P}^k$, and then likewise throughout the construction. But it is an easy exercise to check that the constrained principal angles between $P^k$ and $\overline{P}^k$ coincide with the ordinary principal angles between these complex subspaces of $C^{2n}$. We leave this to the reader, henceforth drop the adjective "constrained", and use this information in what follows.

## (3) PRINCIPAL ANGLES THEOREM FOR CONJUGATE COMPLEX SUBSPACES OF $C^{2n}$.

Let $P^k$ and $\overline{P}^k$ be a pair of conjugate complex subspaces of $C^{2n}$ which meet only at the origin, and $Q^k$ and $\overline{Q}^k$ another such pair. Then there is an element of $O(2n)$ taking $P^k$ to $Q^k$ (and automatically taking $\overline{P}^k$ to $\overline{Q}^k$) if and only if the principal angles between $P^k$ and $\overline{P}^k$ coincide with the principal angles between $Q^k$ and $\overline{Q}^k$.

**Proof.**

Let $P^k$ and $\overline{P}^k$ be a pair of conjugate complex subspaces of $C^{2n}$ which meet only at the origin, and $Q^k$ and $\overline{Q}^k$ another such pair. The condition of matching principal angles is clearly necessary for the existence of an element of $O(2n)$ taking $P^k$ to $Q^k$ and $\overline{P}^k$ to $\overline{Q}^k$.

Suppose, conversely, that the principal angles between $P^k$ and $\overline{P}^k$ coincide with the principal angles between $Q^k$ and $\overline{Q}^k$.

Then by Theorem (2) there is an element $F$ of $U(2n)$ which takes the orthonormal bases

$$v_1, v_2, \ldots, v_k \quad \text{and} \quad \overline{v}_1, \overline{v}_2, \ldots, \overline{v}_k$$

for $P^k$ and $\overline{P}^k$ to the orthonormal bases

$$w_1, w_2, \ldots, w_k \quad \text{and} \quad \overline{w}_1, \overline{w}_2, \ldots, \overline{w}_k$$

for $Q^k$ and $\overline{Q}^k$.

We claim that $F$ commutes with complex conjugation, and hence takes real points of $C^{2n}$ to real points of $C^{2n}$.

Any unit vector in $C\,v_1$ can be written as $e^{i\theta}\,v_1$, and since $F$ is complex linear, $F(e^{i\theta}\,v_1) = e^{i\theta}\,w_1$. Likewise, $F(e^{i\theta}\,\overline{v}_1) = e^{i\theta}\,\overline{w}_1$. This last equality is also true with $\theta$ replaced by $-\theta$, hence $F(e^{-i\theta}\,\overline{v}_1) = e^{-i\theta}\,\overline{w}_1$. But $e^{-i\theta}\,\overline{v}_1$ is the complex conjugate of $e^{i\theta}\,v_1$, and $e^{-i\theta}\,\overline{w}_1$ is the complex conjugate of $e^{i\theta}\,w_1$. Thus $F$ commutes with complex conjugation on $C\,v_1$, and it likewise commutes with complex conjugation on $C\,\overline{v}_1$, so it commutes with complex conjugation on $C\,v_1 + C\,\overline{v}_1$. Similarly, it commutes with complex conjugation on $C\,v_r + C\,\overline{v}_r$, and hence on all of $P^k + \overline{P}^k \to Q^k + \overline{Q}^k$.

If $k = n$, then $P^k + \overline{P}^k$ is all of $C^{2n}$ and so $F$ commutes with complex conjugation on all of $C^{2n}$. If $k < n$, then we can easily modify $F$ on the orthogonal complement of $P^k + \overline{P}^k$ so that it commutes with complex conjugation there as well.

Finally, since $F$ commutes with complex conjugation on all of $C^{2n}$, it takes the real points $R^{2n}$ of $C^{2n}$ to themselves, and is hence an element of the subgroup $O(2n)$ of $U(2n)$.

This completes the proof of (3) above.

## Proof of Proposition 3.

We will exhibit an O(2n)-equivariant deformation retraction of the space of linear complex structures on $R^{2n}$ to its subspace of orthogonal complex structures.

We start with a linear complex structure $J : R^{2n} \to R^{2n}$ and the corresponding direct sum decomposition of the complexification $C^{2n} = V_C^+ + V_C^-$ into a pair of conjugate complex subspaces, the $+i$ and $-i$ eigenspaces of $J : C^{2n} \to C^{2n}$.

We want to move $V_C^+$ and $V_C^-$ apart until they are orthogonal, keeping the intermediate positions as complex conjugates of one another, so as to deform the linear complex structure $J$ through other linear complex structures, until we arrive at the orthogonal complex structure corresponding to the terminal positions of $V_C^+$ and $V_C^-$ in this deformation, as shown below in $C^4$.

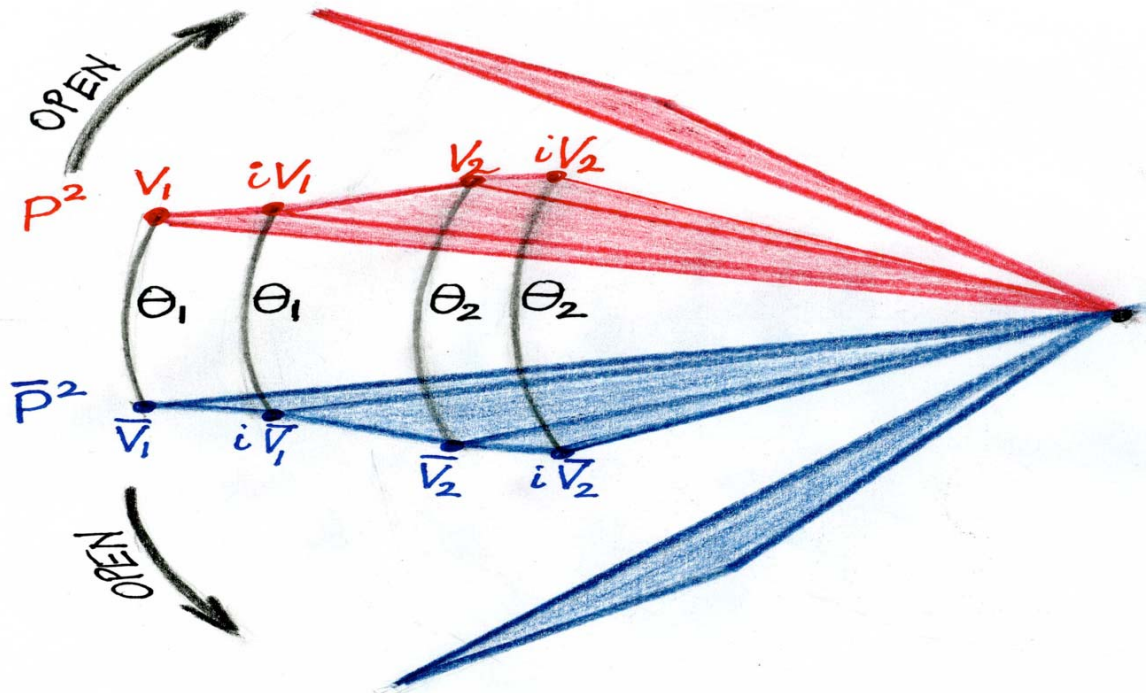

**Figure 10.  Opening up a pair of complex 2-dimensional conjugate subspaces in $C^4$, guided by the principal angles, until they become orthogonal.**

In $C^{2n}$, we open up $V_C^+$ and $V_C^-$ like 2n pairs of scissors in the 2-planes spanned by

$$v_1 \text{ and } \bar{v}_1 \ , \ i\,v_1 \text{ and } i\,\bar{v}_1 \,, \dots , \ v_n \text{ and } \bar{v}_n \ , \ i\,v_n \text{ and } i\,\bar{v}_n \ ,$$

at rates proportional to the complementary angles $\pi/2 - \theta_i$, so that they all open up to angle $\pi/2$ at the same time.

Each of these 2-planes contains a line of real vectors and an orthogonal line of purely imaginary vectors.

As the 2n pairs of scissors open up, the opening vectors $v_k$ and $\bar{v}_k$ remain symmetric with respect to reflection in the real line in their 2-plane, and hence remain conjugates of one another.

By contrast, the opening vectors $i\,v_k$ and $i\,\bar{v}_k$ remain symmetric with respect to reflection in the purely imaginary line in their 2-plane, and hence remain *negative* conjugates of one another.

It follows that the complex 2n-dimensional subspaces $V_C^+$ and $V_C^-$ remain complex conjugates of one another as they open up, until they are finally orthogonal to one another.

This opening up of $V_C^+$ and $V_C^-$ is not affected by the ambiguity in the choice of the above bases for these subspaces, even if several successive principal angles are equal.

During this opening, all the complex structures on $C^{2n}$ commute with complex conjugation, and hence take the subspace $R^{2n}$ of real points to itself.

The result is a deformation retraction of the space of linear complex structures on $R^{2n}$ to its subspace of orthogonal complex structures, and the geometric naturality of all the constructions testifies to the O(2n)-equivariance of this procedure.

This completes the proof of Proposition 3.

# PROOF OF PROPOSITION 4

**PROPOSITION 4.** **There exists a smooth fibration $F$ of $S^{2n+1}$ by oriented great circles whose base space $M_F$ is tangent at $P$ to any preassigned 2n-plane transverse to the bad cone $BC(P)$.**

We begin with a sketch of the proof.

Start in the tangent space $\mathrm{Hom}(P, P^{\perp})$ to $G_2 R^{2n+2}$ at $P$ with a given 2n-plane which is transverse to the bad cone $BC(P)$, hence the graph of a linear map $A: P^{\perp} \to P^{\perp}$ with no real eigenvalues.

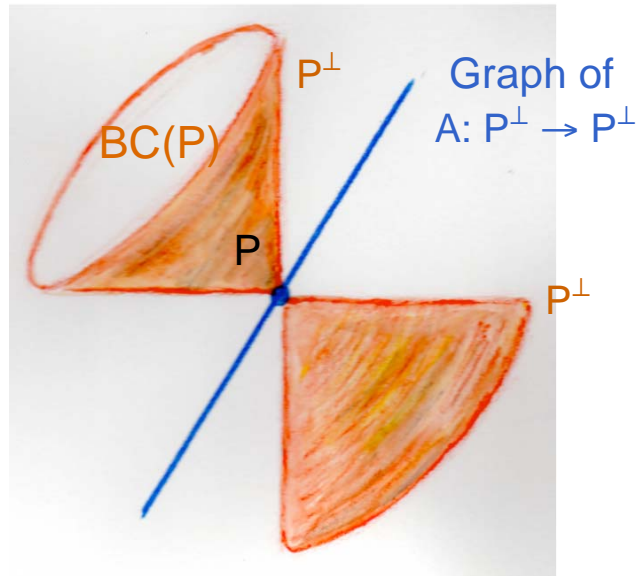

**Figure 11.** **The graph of $A: P^{\perp} \to P^{\perp}$ is transverse to the bad cone**

We must find a fibration $F$ of $S^{2n+1}$ by great circles including $P$, with this preassigned tangent 2n-plane to its base space $M_F$ at $P$.

To do this, let $J_A: P^{\perp} \to P^{\perp}$ be the linear complex structure with the same generalized eigenspaces as $A$, the one to which we deformed $A$ in Proposition 2.

Extend $J_A$ to a complex structure on $R^{2n+2} = P + P^{\perp}$ which rotates the oriented 2-plane $P$ within itself by $90^{o}$.

This complex structure $J_A$ on $R^{2n+2}$ determines a Hopf-like fibration $H_{J_A}$ of $S^{2n+1}$ by the oriented unit circles on the $J_A$-complex lines.

The graph of $J_A: P^\perp \to P^\perp$ is a 2n-plane in $Hom(P, P^\perp)$ which can be regarded as part of the base space of this fibration $H_{J_A}$, and also as its tangent space at $P$.

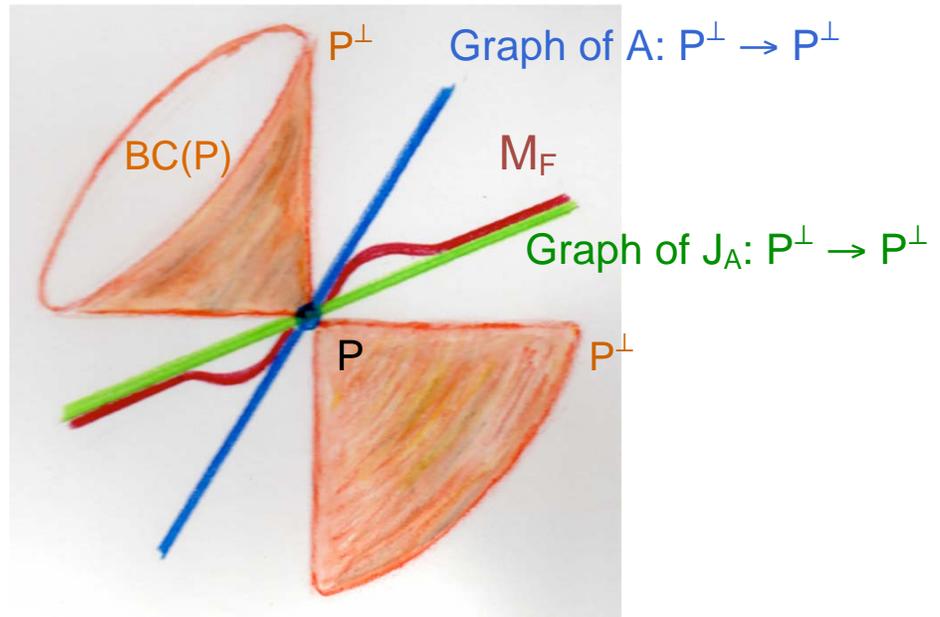

**Figure 12. Interpolating between the graph of A and the graph of the corresponding linear complex structure $J_A$**

We will interpolate between the graphs of $A$ and $J_A$, using the fact that they have the same generalized eigenspaces, to construct the base space $M_F$ of a fibration F of $S^{2n+1}$ by great circles which is tangent at $P$ to the graph of $A$, and which agrees with the fibration $H_{J_A}$ outside a small neighborhood of $P$.

The details of the interpolation are given in the full proof, which we begin now.

## Proof of Proposition 4.

Recall that the 4n-dimensional vector space $\text{Hom}(P, P^\perp)$ serves both as a coordinate neighborhood about $P$ in $G_2 R^{2n+2}$, and as the tangent space to this Grassmannian at $P$.

We start with a 2n-dimensional subspace of $\text{Hom}(P, P^\perp)$ which is the graph of a linear transformation $A: P^\perp \to P^\perp$ with no real eigenvalues. Our goal is to construct a smooth fibration $F$ of $S^{2n+1}$ by oriented great circles, whose base space $M_F$ can be viewed within this neighborhood as the graph of the smooth nonlinear function $N: P^\perp \to P^\perp$, defined by

$$N(x) = f(|x|) A(x) + \big(1 - f(|x|)\big) J(x),$$

for all $x \in P^\perp$; see Figure 12. Here, $f: [0, \infty) \to [0, 1]$ is a smooth bump function which will be defined shortly, and $J = J_A$ is the linear complex structure corresponding to $A$ which was defined in the proof of Proposition 2.

Our task is to choose $f$ so that the differential $dN_x$ of $N$ at each point $x \in P^\perp$ has no real eigenvalues.

We compute $dN_x$ applied to a vector $v$ in $P^\perp$, keeping in mind that the linear functions $A$ and $J$ serve as their own differentials at all points $x$.

$$dN_x(v) = f(|x|) A(v) + \big(1 - f(|x|)\big) J(v)$$
$$+ f'(|x|) (x/|x| \bullet v) A(x) - f'(|x|) (x/|x| \bullet v) J(x).$$

Suppose that $dN_x(v) = \lambda v$ at some point $x \in P^\perp$, for some unit vector $v$, and for some real number $\lambda$.

We will insert this into the previous equation, and then choose the bump function $f$ to prevent this from happening at any point $x$ and for any $\lambda$.

We get

$$\lambda\, v \;=\; f(|x|)\, A(v) \;+\; \big(1 - f(|x|)\big)\, J(v)$$

$$+\; f'(|x|)\, (x/|x| \bullet v)\, A(x) \;-\; f'(|x|)\, (x/|x| \bullet v)\, J(x)\,,$$

and rewrite this as

$$(*) \qquad \lambda\, v \;-\; \big[\, f(|x|)\, A \;+\; \big(1 - f(|x|)\big)\, J\,\big]\, (v)$$

$$=\; f'(|x|)\, (x/|x| \bullet v)\, \big[\, A(x) \;-\; J(x)\,\big]\,.$$

Next we will find an $\varepsilon > 0$ so that the left hand side of $(*)$ has norm $\geq \varepsilon$, independent of the bump function $f$ and the point $x \in P^\perp$. Then we will choose $f$ so that the right hand side has norm $< \varepsilon$.

Suppose first that we cannot find a positive lower bound for the norm of the left hand side.

The left hand side cannot be zero at any $x \in P^\perp$, since the linear maps $t\, A + (1 - t)\, J$ from $P^\perp$ to $P^\perp$ have no real eigenvalues for $0 \leq t \leq 1$, as we showed in the proof of Proposition 2.

Now suppose that as we vary $x \in P^\perp$ among those $x$ for which $dN_x$ has a real eigenvalue, the norm of the left hand side of $(*)$ becomes arbitrarily close to zero. Note that as we vary $x$, the eigenvalue $\lambda$ of $dN_x$, if it exists, might change.

So we suppose that for each integer $n$ there is a real number $\lambda_n$, a unit vector $v_n$ and a real number $t_n \in [0, 1]$ such that

$$|\, \lambda_n\, v_n \;-\; [\, t_n\, A \;+\; (1 - t_n)\, J\,]\, (v_n)\,| \;<\; 1/n\,.$$

We note that the real numbers $\lambda_n$ are bounded in size, since

$$|\, t\, A \;+\; (1 - t)\, J\,| \;\leq\; |\, A\,| \;+\; |\, J\,|$$

is bounded and since $v_n$ is a unit vector.

Then, due to compactness of this bounded interval of real numbers, compactness of the unit 3-sphere in $P^\perp$, and compactness of the interval $[0, 1]$, there is a subsequence $(n_k)$ of the integers with

$$\lambda_{n_k} \to \lambda, \; v_{n_k} \to v \; \text{ and } \; t_{n_k} \to t,$$

so that in the limit we have

$$\lambda v - [\, t A + (1 - t) J\,] (v) = 0,$$

which contradicts the fact that $t A + (1 - t) J$ has no real eigenvalues.

Thus, independent of our choice of $f$ (yet to be made), there is an $\varepsilon > 0$ so that

$$\big| \lambda v - \big[ f(|x|) A + \big(1 - f(|x|)\big) J \big] (v) \big| \geq \varepsilon.$$

We fix this $\varepsilon > 0$ and consider the right hand side of (*),

$$f'(|x|) \, (x/|x| \bullet v) \, \big[ A(x) - J(x) \big],$$

which has norm $\leq |f'(|x|)| \; |A - J| \; |x|$.

We will determine how to choose $f$ so that

$$|f'(s)| \; s < \varepsilon / |A - J|,$$

for any real number $s$ in $[0, \infty)$.

Let $S(f) = \operatorname{Sup} \{\, s \, f'(s) : s \geq 0 \,\}$. We want to choose the bump function $f$ so that $S(f) < \varepsilon / |A - J|$, thus making $S(f)$ as small as necessary.

Start by choosing any smooth bump function $f : [0, \infty) \to [0, 1]$ so that $f(s) = 1$ for $s$ near $0$ and $f(s) = 0$ for $s$ sufficiently large.

Then define $f_n(s) = f(s^{1/n})$ for $n = 1, 2, 3, \ldots$ .

A quick check shows that $S(f_n) = S(f) / n$ , hence for sufficiently large $n$ , the bump function $f_n$ can be used in place of $f$ , so that the right hand side of (*) has norm $< \varepsilon$ .

This shows that (*) is impossible, because the left hand side has norm $\geq \varepsilon$ independent of our choice of $f$ , while for some $f$ , the right hand side has norm $< \varepsilon$ .

This contradicts our supposition that $dN_x(v) = \lambda v$ at some point $x \in P^\perp$ , for some unit vector $v$ , and for some real number $\lambda$ , and therefore confirms that the differential $dN_x$ of $N$ at each point $x \in P^\perp$ has no real eigenvalues.

We now want to define the fibration $F$ of $S^{2n+1}$ by oriented great circles so that its base space $M_F$ within the coordinate neighborhood $\text{Hom}(P, P^\perp)$ is the graph of $N$ , and outside that neighborhood coincides with the base space $M_J$ of the fibration of $S^{2n+1}$ by the unit circles on the J-complex lines.

Since the differential $dN_x$ at each $x \in P^\perp$ has no real eigenvalues, the base space $M_F$ is everywhere transverse to the field of bad cones, and so by Proposition 1 is indeed the base space of a smooth fibration $F$ of $S^{2n+1}$ by oriented great circles.

By construction, $M_F$ agrees with the graph of $A$ near the fibre $P$ , so that we certainly have $T_P M_F = A$ , as required.

This completes the proof of Proposition 4.

# PROOF OF THEOREM A

**THEOREM A.** **The space $\{T_P M_F\}$ of tangent 2n-planes at $P$ to the base spaces $M_F$ of smooth oriented great circle fibrations $F$ of $S^{2n+1}$ deformation retracts to its subspace $\{T_P M_H\}$ of tangent 2n-planes to Hopf fibrations $H$ of $S^{2n+1}$.**

That is, the set of 2n-planes in $T_P G_2 R^{2n+2}$ tangent to the base space of a fibration of $S^{2n+1}$ by great circles deformation retracts to its subspace of 2n-planes tangent to Hopf fibrations.

**Proof.**

Start with the space $\{T_P M_F\}$ of tangent 2n-planes at $P$ to the base spaces $M_F$ of all smooth great circle fibrations $F$ of $S^{2n+1}$.

Use Propositions 1 and 4 to write

$$\{T_P M_F\} = \{\text{ 2n-planes in } T_P(G_2 R^{2n+2}) \text{ transverse to } BC(P) \}$$

$$= \{\text{ Linear maps } T: R^{2n} \to R^{2n} \text{ with no real eigenvalues }\},$$

with $P^\perp$ playing the role of $R^{2n}$.

Then by Propositions 2 and 3, the above space deformation retracts to its subspace

$$\{\text{ Orthogonal complex structures } J: R^{2n} \to R^{2n} \},$$

which is in one-to-one correspondence with the space $\{T_P M_H\}$ of tangent 2n-planes at $P$ to the base spaces $M_H$ of Hopf fibrations $H$ of $S^{2n+1}$ containing the fibre $P$.

This proves Theorem A.

# PROOF OF THEOREM B

**THEOREM B.** **Every germ of a smooth fibration of $S^{2n+1}$ by oriented great circles extends to such a fibration of all of $S^{2n+1}$.**

**Proof.** Let $F$ be a germ of a smooth fibration of $S^{2n+1}$ by great circles containing the fibre P, and $M_F \subset G_2 R^{2n+2}$ its base space.

We must produce a smooth fibration $F''$ of all of $S^{2n+1}$ by great circles which agrees with $F$ in a neighborhood of $P$.

Let $T_P M_F$ be the tangent 2n-plane to $M_F$ at $P$.

We know that $T_P M_F$ is transverse to $BC(P)$, so by Proposition 4, there is a smooth fibration $F'$ of all of $S^{2n+1}$ by great circles with $T_P M_{F'} = T_P M_F$.

By routine interpolation, we get a smooth submanifold $M''$ of $G_2 R^{2n+2}$ which agrees with $M_F$ in a small neighborhood of P, and then agrees with $M_{F'}$ outside a slightly larger neighborhood of $P$, and whose tangent planes are all as close as desired to $T_P M_{F'} = T_P M_F$. See Figure 13.

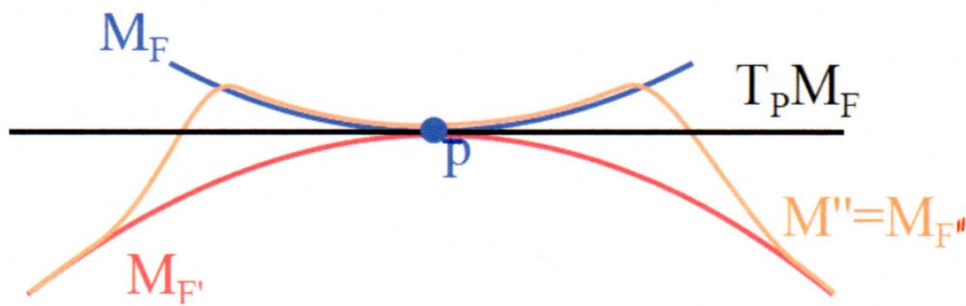

**Figure 13. Interpolation between the base space $M_F$ of the germ and the base space $M_{F'}$ of an entire fibration which is tangent to the germ**

Thanks to this closeness, the tangent planes to $M''$ are transverse to the bad cones at all points, and hence $M'' = M_{F''}$ is the base space of a fibration $F''$ of **all** of $S^{2n+1}$ by great circles. This fibration $F''$ agrees with $F$ in a neighborhood of $P$, completing the proof of Theorem B.

Patricia Cahn
University of Pennsylvania
*pcahn@math.upenn.edu*

Herman Gluck
University of Pennsylvania
*gluck@math.upenn.edu*

Haggai Nuchi
University of Toronto
*hnuchi@math.toronto.edu*